\documentclass[11pt]{article} 
\usepackage{amsthm}
\usepackage{amsmath}
\usepackage{amscd, amssymb}
\usepackage{hyperref}
\usepackage{comment}

\font\rmsml=cmr7

\newtheorem{theorem}{Theorem}[section]
\newtheorem{lemma}[theorem]{Lemma}
\newtheorem{prop}[theorem]{Proposition}
\newtheorem{corollary}[theorem]{Corollary}
\newtheorem{conj}[theorem]{{Conjecture}}

\newtheorem{claim}[theorem]{{Claim}}
\newtheorem{definition}[theorem]{{Definition}}
\newtheorem{hypothesis}[theorem]{{Hypothesis}}
\newtheorem{notation}[theorem]{{Notation}}

\theoremstyle{remark}
\newtheorem{remark}[theorem]{Remark}

\def\bclaim{\begin{claim}}
\def\eclaim{\end{claim}}
\def\bdefin{\begin{definition}}
\def\edefin{\end{definition}}
\def\bcor{\begin{corollary}}
\def\ecor{\end{corollary}}
\def\bthm{\begin{theorem}}
\def\ethm{\end{theorem}}
\def\bconj{\begin{conjecture}}
\def\econj{\end{conjecture}}
\def\blem{\begin{lemma}}
\def\elem{\end{lemma}}
\def\blemma{\begin{lemma}}
\def\elemma{\end{lemma}}
\def\bprop{\begin{prop}}
\def\eprop{\end{prop}}
\def\bremark{\begin{remark}}
\def\eremark{\end{remark}}
\def\bhyp{\begin{hypothesis}}
\def\ehyp{\end{hypothesis}}
\def\bnot{\begin{notation}}
\def\enot{\end{notation}}
\def\bpf{\begin{proof}}
\def\epf{\end{proof}}
\def\beq{\begin{equation}}
\def\eeq{\end{equation}}
\def\beqno{\begin{equation*}}
\def\eeqno{\end{equation*}}
\def\eaeq{\end{aligned}}
\def\baeq{\begin{aligned}}

\def\Aut{\h{\rm Aut}}

\def\JJJ{\h{\rm J}}

\newcommand{\RR}{\mathbb R}

\newcommand{\NN}{\mathbb N}

\renewcommand{\H}{\mathcal H}

\newcommand{\E}{\mathcal E}

\newcommand{\ra}{\right\rangle}

\def\Aut{{\operatorname{Aut}}}

\def\max{{\operatorname{max}}}

\def\q{\quad}
\def\qq{\qquad}

\def\del{\partial}

\def\O{X}

\def\vp{\varphi}

\def\h#1{\hbox{#1}}

\def\lb{\label}

\def\K{K\"ahler }
\def\i{\sqrt{-1}}
\def\del{\partial}
\def\dbar{\bar\partial}
\def\ddbar{\del\dbar}
\def\ra{\rightarrow}

\def\del{\partial}
\newcommand{\calH}{\mathcal{H}}

\def\O{\Omega}

\def\w{\wedge}
\def\o{\omega}

\def\vp{\varphi}
\def\beq{\begin{equation}}
\def\eeq{\end{equation}}
\def\bi#1{\bibitem{#1}}
\def\PSH{\mathrm{PSH}}

\def\h#1{\hbox{#1}}

\def\calE{{\mathcal E}}

\def\K{K\"ahler } \def\Kno{K\"ahler}
\def\KE{K\"ahler--Einstein } \def\KEno{K\"ahler--Einstein}

\def\polishl{\char'40l}

\def\psiKE{\psi_{\h{\rmsml KE}}}

\def\ra{\rightarrow}
\def\pa{\partial}

\def\w{\wedge}
\def\i{\sqrt{-1}}
\def\Aut{{\operatorname{Aut}}}

\def\Cinf{C^\infty}

\def\del{\partial}

\def\ovp{{\o_{\vp}}}
\def\on{\omega^n}
\def\ovpn{\omega_{\vp}^n}

\def\intM{\int_M}
\def\Ric{\hbox{\rm Ric}\,}

\def\intM{\int_M}  \def\intm{\int_M}
\def\V{\frac1V}

\def\ovp{\o_{\vp}}
\def\on{\o^n}
\def\ovpn{\o_{\vp}^n}
\def\Ric{\h{\rm Ric}}

\def\Ricovp{\Ric\,\ovp}

\def\Cinf{\h{\cal C}^\infty}
\def\CinfM{\Cinf(M)}

\def\Ho{{\calH_{\o}}}
\def\AM{\operatorname{AM}}

\title{Convergence of the K\"ahler--Ricci iteration}
\author{Tam\'as Darvas and Yanir A. Rubinstein}
\date{}

\begin{document}
\maketitle

\begin{abstract}
The Ricci iteration is a discrete analogue of the Ricci flow.
According to Perelman, the Ricci flow converges to 
a K\"ahler--Einstein metric whenever one exists, 
and it has been conjectured that the Ricci iteration should behave similarly. This article confirms this conjecture. 
As a special case, this gives a new method of uniformization of the Riemann sphere.
\end{abstract}

\def\lb{\label}

\section{Introduction}

Let $(M,g_1)$ be a compact Riemannian manifold. 
A Ricci iteration is a sequence of metrics $\{g_i\}_{i\in\mathbb N}$ on $M$ satisfying 
\begin{align}
\label{RIEq}
\Ric g_{i+1}=g_i,\q i\in\mathbb N,
\end{align}
where $\Ric g_{i+1}$ denotes the Ricci curvature of $g_{i+1}$. 
One may think of~\eqref{RIEq}
as a dynamical system on the space of Riemannian  metrics on $M$. 
Part of the interest in the Ricci iteration is that, clearly, 
Einstein metrics are fixed points, and so~\eqref{RIEq} aims to provide a 
natural theoretical and numerical 
approach to uniformization in the challenging case of positive
Ricci curvature (different Ricci iterations can be defined in the context of non-positive curvature, but these are typically easier to understand
and will not be discussed here).
In essence, the Ricci iteration aims to reduce the Einstein equation to a sequence of prescribed Ricci curvature equations and can be thought of as a  discretization of the Ricci flow.
Going back to \cite{R07,R08}, it has been studied since by a number of 
authors \cite{Berm,BBEGZ,CPS,CS,CS2,CW,GKY,JMR,Ke,PR}, 
see also the survey~\cite[\S6.5]{R14}.

Of particular interest has been the study of the Ricci iteration on K\"ahler manifolds
(for the non-\K case results are scarce, see \cite{PR}).
When $(M,\JJJ,g_1)$ is K\"ahler, the Calabi--Yau Theorem
\cite{Yau1978} guarantees the existence and uniqueness of the sequence $\{g_i\}_{i\in\NN}$ if and only if  $M$ is Fano (i.e., has positive first Chern class $c_1(M,\JJJ)$)
and the \K class associated to $g_1$ is $c_1(M,\JJJ)$. Under a rather restrictive technical assumption, 
one of us showed that $g_i$ converges smoothly to a K\"ahler--Einstein metric 
\cite[Theorem 3.3]{R08} and made the following general conjecture
\cite[Conjecture 3.2]{R08}:  

\begin{conj}
\lb{RConj}
Let $(M,\JJJ,g_1)$ be a compact \K manifold 
admitting a \KE metric. Suppose the \K class associated to $g_1$ is $c_1(M,\JJJ)$.
Then  the  Ricci iteration \eqref{RIEq} converges in the sense of Cheeger-Gromov 
to a \KE metric.
\end{conj}

The best result so far on this conjecture is due to Berman et al.
\cite{BBEGZ} who replace the technical assumption of \cite[Theorem 3.3]{R08}
concerning Tian's $\alpha$-invariant by the weaker assumption of the Mabuchi energy being proper (both of these assumptions imply a \KE metric exists). 
Therefore, by a classical result of Tian \cite{T97}, Conjecture \ref{RConj} holds if 
$M$ admits no holomorphic vector fields. However, the properness assumption is still too restrictive and fails  in general. For example,
Conjecture \ref{RConj} is still open even for $M=S^2$, the two-sphere. 
Furthermore, as recent counterexamples show \cite{DR}, it is not possible to modify the properness assumption to simply hold on $K$-invariant metrics, where $K$ is the maximal compact subgroup of the holomorphic automorphism group of $M$.

The main result of the present article is the resolution of Conjecture \ref{RConj},
and in fact with a stronger convergence.

\begin{theorem}
\lb{MainThm}
Let $(M,\JJJ,g_1)$ be a compact \K manifold 
admitting a \KE metric. Suppose the \K class associated to $g_1$ is $c_1(M,\JJJ)$
and let $\{g_i\}_{i\in\NN}$ be given by \eqref{RIEq}. 
Then  there exists holomorphic diffeomorphisms $h_k$ such that 
$h_k^*g_k$ converges smoothly to a \KE metric.
\end{theorem}

\font\small=cmr5                       

\def\oFS{\o_{\hbox{\small FS}}}
\def\oFSc{\o_{\hbox{\small FS},c}}
\def\dbz{d\bar z}

\subsection{Uniformization of the two-sphere}
\label{}

As a very special case we obtain the following new method of uniformization. 
Fix a conformal class of volume $V$ on $S^2$.
As we know, in this class there is a constant curvature metric, the
round one.
More precisely, let $\o_c$ denote the round
form of the constant $c$ Ricci curvature metric on $M=(S^2,\JJJ)$, given locally by
$$
\o_c=\frac\i{c\pi}\frac{dz\w\dbz}{(1+|z|^2)^2}.
$$
Here $V=\int_{S^2}\o_c=c_1([M])/c=2/c$. Consequently, 
$c=1/2\pi$ in case we are restricting the Euclidean metric of $\RR^3$ to
the unit sphere.

Let $\omega$ be any metric on $S^2$ with $\int_{S^2} \omega = V = 2/c$. Introduce $u_0=0$, and we solve iteratively to find
$u_i\in C^\infty(S^2)$ satisfying
\beq
\lb{RI2DEq}
\Delta_{\omega} u_i = R_\omega-2e^{u_{i-1}},
\ \ \ \ 
\h{and\ \ \  $\int_{S^2}e^{u_i}\omega=2/c$},
\eeq
so that the scalar curvature of $\o_i:=e^{u_i}\omega$
satisfies $R_{\o_i}=2 e^{u_{i-1}-u_i}$, or equivalently, $\Ric\o_i=\o_{i-1}$.
(In two dimensions, $\Ric\,\omega = \frac12 R_\omega \, \omega$, where $R_\omega$ is the scalar curvature. If $\omega_1 = e^{\phi}\omega_0$, then the scalar curvatures of these two metrics satisfy 
$$
\Delta_{\omega_0} \phi -  R_{\omega_0} +  R_{\omega_1} e^{\phi} = 0.
$$
We note that the conformal factor is
often written $e^{2\phi}$ elsewhere, but this is compensated for here by the fact that $R_\omega = 2K_\omega$, where $K_\omega$  is the Gauss curvature.)

\bcor We fix $c>0$ and let $\omega$ be any K\"ahler form on $S^2$ with $\int_X \omega =2/c$. We introduce $\{u_i\}\subset
C^\infty(S^2)$ by repeatedly solving the Poisson equation \eqref{RI2DEq}.
Then, there exist M\"obius transformations $h_i$ such that $h_i^*(e^{u_i}\omega)$ converges 
smoothly to the round metric $\o_c$.
\ecor

\subsection{Discretization of the Ricci flow}
\label{}

One of the original motivations for introducing the Ricci iteration, 
going back to \cite{R07,R08}, is its relation to the Ricci flow.
Hamilton's Ricci flow on a \K manifold of definite or zero first Chern class is defined as $\{\o(t)\}_{t\in\RR_+}$ satisfying the evolution equation
$$
\begin{aligned}
\frac{\pa \o(t)}{\pa t} & =-\Ric\,\o(t)+\mu\o(t),\quad t\in\RR_+,\cr
\o(0) & =\o,
\end{aligned}
$$
where $\O$ is a \K class satisfying $\mu\O=c_1(M,\JJJ)$ for $\mu\in\{-1,0,1 \}$ 
and $[\o]=\O$
\cite{Ham82}.

The following dynamical system is seen to be a discrete version of this flow \cite[Definition 3.1]{R08}, obtained by a backward Euler discretization with time step $\tau$.

\begin{definition}
\label{RIdefinition}
Let $\O$ be a \K class satisfying $\mu\O=c_1(M,\JJJ)$ for $\mu\in\{-1,0,1 \}$.
Given a \K form $\o$ with $[\o]=\O$ and a number $\tau>0$, define 
the time $\tau$ Ricci iteration %
 to be 
the sequence of forms $\{\o_{k\tau}\}_{k\ge0}$ 
satisfying the equations
$$
\begin{aligned}
\frac{\o_{k\tau}-\o_{(k-1)\tau}}{\tau} & =-\Ric\,\o_{k\tau}+\mu\o_{k\tau},\quad   k\in\NN,\cr
\o_0 & =\o.
\end{aligned}
$$
\end{definition}

Let us assume that $\mu=1$ from now on
(for the cases $\mu\in \{-1,0\}$ see 
\cite[Theorem 3.3]{R08}).
Observe that in the case when $\tau=1$, the time $\tau$ Ricci iteration is precisely the Ricci iteration from \eqref{RIEq}.
Indeed, Conjecture \ref{RConj} is 
in fact a special case of the following conjecture concerning
the time $\tau$ Ricci iteration for any $\tau>0$
\cite[Conjecture 3.2]{R08}.

\begin{conj}
\lb{MainConj2}
Let $(M,\JJJ)$ be a compact \K manifold admitting a \KE metric. 
Let $\O$ be a \K class such that $\O=c_1(M,\JJJ)$.
Then for any $\o$ with $[\o]=\O$ and for any $\tau>0$, the time $\tau$ Ricci iteration exists for all $k\in\NN$ and converges
in the sense of Cheeger-Gromov 
to a \KE metric.
\end{conj}

The case when $\tau>1$ is treated in  \cite[Theorem 3.3]{R08}. 
However, it is the case $\tau\le1$ that is the most interesting and challenging.
The case $\tau=1$ is perhaps the most interesting due to the simple
geometrical interpretation \eqref{RIEq} while the cases $\tau<1$ are
interesting due to the connection to the \Kno--Ricci flow. 
 In this regime one may expect the Ricci iteration to converge to the Ricci flow in a certain scaling limit as $\tau\ra0$.
The cases $\tau\le1$
are challenging since the a priori estimates are considerably
harder then. While in the regime $\tau>1$ one has a uniform positive
Ricci lower bound along the iteration, this is no longer true when
$\tau\le 1$. Thus, there is no a priori control on the diameter 
or the Poincar\`e
and Sobolev constants. We work around these difficulties, by analyzing the Ricci iteration in the metric geometry of the space of K\"ahler potentials \cite{Da}.

In this article we in fact confirm the more general Conjecture \ref{MainConj2},
and treat the iteration for all time steps $\tau$
by proving the following result of which 
Theorem \ref{MainThm} is a special case.

\begin{theorem}
\lb{MainThm2}
Let $(M,\JJJ,g_1)$ be a compact \K manifold 
admitting a \KE metric. Suppose the \K class associated to $g_1$ is $c_1(M,\JJJ)$
and let $\{\o_{k\tau}\}_{k\in\NN}$ be the time $\tau$ Ricci iteration
given by Definition \ref{RIdefinition}. 
Then  there exists holomorphic diffeomorphisms $h_k$ such that 
$h_k^*\o_{k\tau}$ converges smoothly to a \KE form.
\end{theorem}

\section{Energy functionals} 
\lb{EnergySec}

Let $(M,\o)$ denote a connected compact closed \K manifold. 
The space of smooth strictly $\o$-plurisubharmonic 
functions (\K potentials)
\beq
\label{HoEq}
\calH_\o 
:=
\{\vp\in C^{\infty}(M) \,:\, \,  \omega_\vp:=\o+\i\ddbar\vp>0\},
\end{equation} 
can be identified with $\H\times\RR$, where
\beq
\label{HEq}
\calH
=
\{\omega_\vp \,:\, \vp\in C^{\infty}(M), \,  \omega_\vp>0\}
\end{equation}
is the space of all \K metrics (or forms) representing the fixed cohomology class $[\o]$.

From now on let $\o$ be a K\"ahler form on $M$, cohomologous to $c_1(M,\JJJ)$. The Aubin--Mabuchi functional was introduced by Mabuchi \cite[Theorem 2.3]{Mabuchi},
\begin{equation}
\label{AMdef}
\h{\rm AM}(\vp):=
\frac{V^{-1}}{n+1}\sum_{j=0}^{n}\int_M \vp\, \o^j \wedge \o_\vp^{n-j},
\end{equation}
where  $V:=\int_M\o_\vp^n=\int_M\o_\vp^n$ is the total volume of the K\"ahler class. Integration by parts gives the useful estimates
\begin{equation}\label{eq:AMest}
\frac{1}{V}\int_M (u - v)\o_u^n \leq \AM(u) - \AM(v)
 \leq \frac{1}{V}\int_M (u - v)\o_v^n. \end{equation}
The subspace
\beq
\lb{H0Eq}
\calH_0:=\h{\rm AM}^{-1}(0)\cap \Ho
\eeq
is isomorphic to $\calH$ \eqref{HEq}, the space of K\"ahler metrics.


Let $f_{\ovp}\in\CinfM$ denote 
the unique function (called the Ricci potential of $\ovp$) satisfying
$$
\i\ddbar f_{\ovp}=\Ricovp-\ovp,\q
\V\intm e^{f_{\ovp}}\ovpn=1.
$$
The Ding and Mabuchi functionals are given by \cite{Ding,Mabuchi}
\beq
\baeq
\lb{FEEq}
D(\vp)
&:=
-\AM(\vp)-\log\V\intM e^{f_\o-\vp}\o^n,
\cr
E(\vp)
&:= \frac{1}{V}\int_X \log \frac{\o_\vp^n}{e^{ f_\o}\o^n} \o_\vp^n - \AM(\vp) + \frac{1}{V}\int_M \vp \o_\vp^n + \frac{1}{V}\int_M f_\o \omega^n. \eaeq
\eeq
Notice that these functionals are invariant
under addition of constants to $\vp$, hence they descend to $\mathcal H$. Additionally, the critical points of these functionals are exactly the \KE metrics.

For $\vp \in \mathcal H_\o$ with $\intM e^{f_\o-\vp}\o^n=V$, 
Jensen's inequality for the convex weight $t \to t \log t$ yields,
\begin{equation}\label{eq: Jensen_ineq}
\textup{Ent} (e^{f_\o-\vp}\o^n,\omega^n_\vp):=
\frac{1}{V}\int_X \log \frac{\o_\vp^n}{e^{f_\o-\vp}\o^n} \o_\vp^n  
=
\frac{1}{V}\int_X 
\frac{\o_\vp^n}{e^{f_\o-\vp}\o^n}
\log \frac{\o_\vp^n}{e^{f_\o-\vp}\o^n}
e^{f_\o-\vp}\o^n  \geq 0.
\end{equation}
Thus,
$$
E(\ovp)
-\V\intM f_\o \on = \textup{Ent} (e^{f_\o-\vp}\o^n,\omega^n_\vp)-\AM(\vp) \geq -\AM(\vp)=D(\ovp).
$$
Moreover, if 
$$
D(\ovp)=
E(\ovp)
-\V\intM f_\o \on
$$
then equality holds in \eqref{eq: Jensen_ineq}. As a result, $\o_\vp^n = e^{f_\o-\vp}\o^n = e^{f_{\o_\vp}}\o_\vp^n$, i.e., 
$\ovp$ is \KEno.
This together with the fact that \KE metrics minimize both $D$ and $E$ allows to conclude the following result (see also \cite[(24)]{R08JFA}):
\bprop 
\lb{RProp}
For $\vp\in\H_\omega$,
$$
D(\ovp)\le E(\ovp)
-\V\intM f_\o \on,
$$
with equality if and only if $\Ricovp=\ovp$. 
%
\eprop

\section{The metric completion}
\label{}

All of the functionals introduced in the previous section can be extended to the potential space $\mathcal E_1$ introduced by Guedj--Zeriahi \cite{GZ}, that can be identified with a natural metric completion of $\calH$ \cite{Da}. 
The resulting metric theory provides essential tools for proving our main result concerning convergence of the Ricci iteration.
We briefly recall this machinery, referring to \cite[\S4--5]{DR} and references therein for more details.

Let
$$
\PSH(M,\o)=\{ \vp \in L^1(M,\o^n)\,:\, \h{$\vp$ is upper semicontinuous and } 
\ovp\geq 0 \}.
$$
Following Guedj--Zeriahi \cite[Definition 1.1]{GZ} we define the subset of full mass potentials:
$$
\E(M,\o):=
\big\{
\vp\in \PSH(M,\o)\,:\,
\lim_{j\ra-\infty}\int_{\{\vp\le j\}}(\o+\i\ddbar\max\{\vp,j\})^n= 0
\big\}.
$$
For each $\vp\in\E(M,\o)$, define
$
\ovpn:=\lim_{j\ra-\infty}
{\bf 1}_{\{\vp>j\}}
(\o+\i\ddbar\max\{\vp,j\})^n.
$
By definition, ${\bf 1}_{\{\vp>j\}}(x)$ is equal to $1$ if $\vp(x)>j$ and zero otherwise,
and
the measure $(\o+\i\ddbar\max\{\vp,j\})^n$ is defined by the
work of Bedford--Taylor \cite{BT} since $\max\{\vp,j\}$ is bounded. Consequently, $\vp \in \mathcal E(M,\omega)$ if and only if $\int_X  \omega_\vp^n = \int_X \omega^n$, justifying the name of $\mathcal E(M,\omega)$. 

Next, define a further subset, the space of finite 1-energy potentials:
$$
\E_1:=\big\{\vp\in\E(M,\omega)\,:\, \int|\vp|\ovpn<\infty\big\}.
$$
Consider the following weak 
Finsler metric on $\mathcal H_\o$ \cite{Da}:
\begin{equation}\label{FinslerDef}
\|\xi\|_\vp:=  V^{-1}\int_M |\xi| \ovpn, 
\q  \xi \in T_\vp \mathcal{H}_\o=C^\infty(M).
\end{equation}
We denote by $d_1$ the associated pseudo-metric and recall the result alluded to above, characterizing the $d_1$-metric completion  of $\mathcal H_\omega$ \cite[Theorem 2, Theorem 3.5]{Da}:
\bthm
\lb{d1CompletionThm}
$(\Ho, d_1)$ is a metric space whose completion can be identified with $(\E_1,{d_1})$,
where 
$$
d_1(u_0,u_1):=\lim_{k\ra\infty}
d_1(u_0(k),u_1(k)),
$$
for any smooth decreasing sequences $\{u_i(k)\}_{k\in\NN}\subset\Ho$
converging pointwise to $u_i \in \mathcal E_1, i=0,1$.
\ethm

Also, by \cite[Theorem 3]{Da}, we have the following qualitative estimates  for the $d_1$ metric in terms of analytic quantities:
\begin{equation}\label{eq: d_1_double_est}
\frac{1}{C}d_1(u,v) \leq \int_M |u-v|\o_u^n + \int_M |u-v|\o_v^n \leq Cd_1(u,v), \ \ u,v \in \mathcal E_1,
\end{equation}
where $C>1$ only depends on $\o$.

A crucial fact is that the formulas defining the energy functionals discussed in 
\S\ref{EnergySec} actually make sense on the metric completion $\mathcal E_1$, and then coincide with the greatest lower semi-continuous extension of the said functionals restricted to $\H_\o$
\cite[Lemma 5.2, Proposition 5.19, Proposition 5.21]{DR}:

\blem
\lb{E1capH0Lemma}
(i)
$\h{\rm AM},D: \mathcal H_\o \to \Bbb R$ each admit a unique $d_1$-continuous extension 
to $\E_1$ 
and these extensions still satisfy \eqref{AMdef} and \eqref{FEEq} respectively.\\
\noindent
(ii)
$E: \mathcal H_\o \to \Bbb R$ admits a $d_1$-lower semi-continuous extension 
to $\E_1$ 
and the greatest such extension still satisfies \eqref{FEEq}.
\elem

Proposition \ref{RProp} was generalized by Berman \cite[Theorem 1.1]{Berm} to the context of the metric completion (for a proof using the Ricci iteration see \cite[Proposition 4.42]{Dasurvey}):

\bthm\lb{BThm}
Proposition \ref{RProp} holds more generally
for all $\vp\in\calE_1$.
\ethm

Let $G:=\Aut_0(M) $ denote the connected component of the complex Lie group of automorphisms (biholomorphisms) of $M$.
The automorphism group  acts on $\H$ by pullback:
\beq\lb{AutActionEq}
f.\eta:=f^\star\eta, \qq f\in G, \q \eta\in\H.
\eeq
Given the one-to-one correspondence between $\mathcal H$ and $\mathcal H_0$
(recall \eqref{H0Eq}), the group $G$  also acts on $\H_0$. The precise action is described in the next lemma \cite[Lemma 5.8]{DR}.
\blem
\lb{AutLem}
For $\vp \in \mathcal H_0$ and $f \in G$ let $f.\vp\in\H_0$ be the unique potential such that $f^*\ovp=\o_{f.\vp}$.
Then,
\beq
\lb{factionAMEq}
f.\vp=f.0+\vp\circ f.
\eeq
\elem
Complementing the above, $G$ acts on $\mathcal H_0$ by $d_1$-isometries \cite[Lemma 5.9]{DR}, which allows to introduce a natural (pseudo)metric on the space $\mathcal H_0/G$:
\begin{equation}\label{eq: d_1_G_def}
d_{1,G}(Gu,Gv)= \inf_{g \in G}d_1(u,g.v), \ \ u,v \in \mathcal H_0.
\end{equation}

\section{Metric convergence of the iteration}
\label{}

We consider the $\tau$-step Ricci iteration equation:
$$\frac{\o_{\psi_{(k+1)\tau}}-\o_{\psi_{k\tau}}}{\tau}=\o_{\psi_{(k+1)\tau}} 
- \Ric\, \o_{\psi_{(k+1)\tau}},$$
for $\tau \in (0,1]$. 
When $\tau =1$, the iteration simply becomes $\Ric\, \o_{\psi_{k+1}} = 
\o_{\psi_k}$. As explained in \cite[(33)]{R08}, on the level of scalars the iteration can be written in the following manner:
\begin{equation}
\label{RIPotEq}
\o_{\psi_{(k+1)\tau}}^n
= 
e^{
f_\o-\frac{1}{\tau}\psi_{k\tau} 
-\big(1- \frac{1}{\tau}\big)\psi_{(k+1)\tau}
}\on,\q k\in\NN,
\end{equation}
with the natural normalization
\beq
\lb{psiknormEq}
\V\intm 
e^{
f_\o-\frac{1}{\tau}\psi_{k\tau} 
-\big(1- \frac{1}{\tau}\big)\psi_{(k+1)\tau}
}
\on=1.
\eeq
Other normalizations may be considered on the level of scalars. In our particular case, there will be special emphasis on working in the geodesically complete potential space $\mathcal H_0$, and we introduce accordingly:
\begin{equation}\label{eq: psi'_k_def}
\psi'_{k\tau} := \psi_{k\tau} - \AM(\psi_{k\tau}) \in \mathcal H_0.
\end{equation}

First we generalize an inequality of \cite{R08} (in the case $\tau =1$) that provides a comparison of the Ding and Mabuchi energies along the $\tau$-iteration:
\begin{prop} \label{prop: tau_energy_ineq}Suppose $\tau \in (0,1]$ and $(M,\o_{\psi^\tau_1})$ is a compact Fano manifold. Then the following estimate holds:
\begin{equation}\label{eq: tau_energy_ineq}
E(\o_{\psi_{(k+1)\tau}}) 
-
\frac{1}{V}\int_M f_{\o}\o^n 
\leq 
\frac{1}{\tau} D(\o_{\psi_{k\tau}}) 
+ 
\Big(1 - \frac{1}{\tau}  \Big)D(\o_{\psi{(k+1)\tau}}), \ \forall \  k \in \Bbb N.
\end{equation}
\end{prop}

In the argument below (and thereafter) we will suppress the parameter $\tau$ from superscripts whenever this will cause no confusion.

\begin{proof} 
Using \eqref{FEEq} and \eqref{RIPotEq},
\begin{flalign*}
E(\o_{\psi_{k+1}})-\frac{1}{V}& \int_M  f_\o \omega^n = \frac{1}{V}\int_X \log \frac{\o_{\psi_{k+1}}^n}{e^{ f_\o}\o^n} \o_{\psi_{k+1}}^n - \AM(\psi_{k+1}) + \frac{1}{V}\int_M {\psi_{k+1}} \o_{\psi_{k+1}}^n\\
&= -\frac{1}{V}\int_M \bigg(\frac{1}{\tau}\psi_{k} +\bigg(1- \frac{1}{\tau}\bigg)\psi_{k+1} \bigg)\o_{\psi_{k+1}}^n - \AM(\psi_{k+1}) + \frac{1}{V}\int_M {\psi_{k+1}} \o_{\psi_{k+1}}^n\\
&=\frac{1}{\tau V}\int_M (\psi_{k+1}-\psi_{k} )\o_{\psi_{k+1}}^n - \AM(\psi_{k+1}).
\end{flalign*}
Using this identity, to finish the proof, we notice that it is enough to prove the following two inequalities (and later add them up):
\begin{equation}\label{eq: AM_ineq}
\frac{1}{\tau V}\int_M (\psi_{k+1}-\psi_{k} )\o_{\psi_{k+1}}^n - \AM(\psi_{k+1}) \leq -\frac{1}{\tau} \AM(\psi_k) - \Big(1-\frac{1}{\tau} \Big) \AM(\psi_{k+1})
\end{equation}
\begin{equation}\label{eq: exp_ineq}
0 \leq -\frac{1}{\tau}\log \bigg(\frac{1}{V}\int_M e^{f_\o - \psi_k}\o^n \bigg) - \Big(1-\frac{1}{\tau}\Big)\log \bigg( \frac{1}{V}\int_M e^{f_\o - \psi_{k+1}}\o^n\bigg)
\end{equation}
Notice that, after rearranging terms, \eqref{eq: AM_ineq} is seen to be equivalent to 
$$
\frac{1}{V} \int_M (\psi_{k+1}-\psi_{k} )\o_{\psi_{k+1}}^n 
\leq \AM(\psi_{k+1}) - \AM(\psi_{k}).
$$ 
Thus,
\eqref{eq: AM_ineq} follows from \eqref{eq:AMest}.
To address \eqref{eq: exp_ineq} 
 we prove the following more general claim.

\begin{claim} For $\tau \in (0,1]$ and $g,h \in C^\infty(X)$ the following estimate holds:
\begin{equation}\label{eq: Holder_ineq}
\bigg(\frac{1}{V}\int_M e^{f_\o - g}\o^n \bigg)^{\frac{1}{\tau}}\bigg( \frac{1}{V}\int_M e^{f_\o - h}\o^n\bigg)^{1 - \frac{1}{\tau}} \leq \frac{1}{V} \int_M e^{f_\o - \frac{1}{\tau}g -\big(1 - \frac{1}{\tau} \big) h}\o^n.
\end{equation}
By our choice of normalization \eqref{psiknormEq}, this  inequality implies \eqref{eq: exp_ineq}.
\end{claim}
As \eqref{eq: Holder_ineq} is seen to be invariant under adding constants to $g$ and $h$, we can assume that $\frac{1}{V}\int_M e^{f_\o - h} \o^n = 1$. In particular, we only have to argue that 
$$\bigg(\frac{1}{V}\int_M e^{-g + h} e^{f_\o - h}\o^n \bigg)^{\frac{1}{\tau}} \leq \frac{1}{V} \int_M \big( e^{-g + h}\big)^{\frac{1}{\tau}} e^{f_\o - h} \o^n.
$$
This follows from Jensen's inequality, as the function $f(t)=t^{\frac{1}{\tau}}$ is convex for $t >0$.
\end{proof}

Next we show that in case a K\"ahler--Einstein metric exists, the iteration $\{\psi'_k\}_{k}$ $d_1$-converges up to pullbacks:
\bprop
\lb{weakProp}
Let $\tau \in (0,1]$.
Suppose a \KE metric exists in $\mathcal H$, and let $\{\psi_{k\tau}\}_{k\in\NN}$ be the solutions of
\eqref{RIPotEq}. Then there exist $g_k \in G$ such that 
$g_k.\psi'_{k\tau}$ $d_1$-converges to a \KE potential.
\eprop

\bpf
Proposition \ref{prop: tau_energy_ineq} combined with Proposition \ref{RProp} gives 
\beq
\lb{MonInEq}
D(\o_{\psi_{k+1}})
\le
E(\o_{\psi_{k+1}}) -\frac{1}{V}\int_M f_{\o}\o^n \leq \frac{1}{\tau} D(\o_{\psi_k}) + \Big(1 - \frac{1}{\tau}  \Big)D(\o_{\psi_{k+1}}), \q k\in\NN.
\eeq
As a result, $\{D(\o_{\psi_{l}})\}_l$ is a decreasing sequence
(this is proved in \cite[Proposition 4.2(ii)]{R08} for $\tau=1$). We fix a \KE potential 
$$
\psiKE \in \mathcal H_0.
$$ 
Existence of such a potential implies that both $D$ and $E$ are
bounded below \cite{BM,DT}. Thus, the (monotone) sequence
$\{D(\o_{\psi_{l}})\}_l$ converges. 
By \eqref{MonInEq}, 
$\{ E(\o_{\psi_{l}})-\frac{1}{V}\int_M f_\o \o^n\}_l$ converges too and both of these sequences have the same limit $l \in \Bbb R$. 

Next we focus on the potentials $\psi'_l \in \mathcal H_0$. By \cite[Theorem 2.4]{DR}, $E$ is $G$-invariant and 
$$
E(\psi_{l}')\ge C_1d_{1,G}(0,\psi_l')-C_2,
$$
and so $d_{1,G}(0,\psi'_l)\le C'$. By definition (see \eqref{eq: d_1_G_def}), there exists $g_l\in G$ such that 
\begin{equation}\label{eq: d_1est_pullback}
d_{1}(\psiKE,g_l.\psi'_l)\le d_{1,G}(G\psiKE,G\psi'_l) + \frac{1}{l} \leq C' + 1.
\end{equation}

\bremark
In fact, there exists $g_l$ which achieve the equality
$d_{1}(\psiKE,g_l.\psi'_l)=
d_{1,G}(G\psiKE,G\psi'_l)$
by \cite[Proposition 6.8]{DR} but we do not have
to know that for our proof here.
\eremark

Denoting 
$$
v_l:= g_l. \psi'_l,
$$ by $G$-invariance of $E$, we obtain that $E(v_l)$ is  bounded. On the other hand, a combination of  \eqref{eq: d_1_double_est} and \eqref{eq: d_1est_pullback} gives that $\AM(v_l)=0$ and $\int_M v_l \omega_{v_l}^n$ are  bounded as well. Comparing with \eqref{HEq}, we see that  $\textup{Ent}(e^{f_0}\o^n,\o_{v_l}^n)$ is bounded too.

By \eqref{eq: d_1_double_est}, $d_1$-boundedness of potentials implies $L^1$-boundedness, which in turn implies boundedness of the supremum.
As a result, we can apply the compactness result of \cite{BBEGZ} 
(see \cite[Theorem 5.6]{DR} for a convenient formulation for our context) 
to conclude that $\{v_l\}_l$ is $d_1$-precompact. 

Next we claim that $d_1(\psiKE,v_l) \to 0$. If this is not the case, then by possibly choosing a subsequence, we can assume that $d_1(\psiKE,v_l)>\varepsilon >0$. 
By possibly choosing another subsequence,  we can assume that $d_1(v_l,u) \to 0$ for some $u \in \mathcal E_1$. Lemma \ref{E1capH0Lemma} gives that $l=D(u)=E(u)-\frac{1}{V}\int_M f_\o \o^n$, in particular $u$ is a \KE potential by Theorem \ref{BThm}. 

By the Bando--Mabuchi uniqueness theorem $u = h.\psiKE$ for some $h \in G$ \cite{BM}. Combining this with \eqref{eq: d_1est_pullback}, we conclude that
$$ d_1(v_{k_l},\psiKE) -\frac{1}{k_l}\leq d_{1,G}(Gv_l,G\psiKE)\leq d_1(h^{-1} v_l,\psiKE)=d_1(v_l,h.\psiKE) = d_1(v_l,u).$$
By choice, the right hand side converges to zero, and the $\liminf$ of left hand side is bounded below by  $\varepsilon >0$, giving a contradiction. This implies that $d_1(v_k,\psiKE) \to 0$, concluding the proof.
\epf

\section{A priori estimates and smooth convergence}

In this section we prove our main result 
by strengthening Proposition \ref{weakProp}.

\bthm
\lb{strongthm}
Let $\tau \in (0,1]$.
Suppose a \KE metric exists in $\mathcal H$, 
and let $\{\psi_{k\tau}\}_{k\in\NN}$ be the solutions of
\eqref{RIPotEq}. Then there exist $g_k \in G$ such that 
$g_k.\psi'_{k\tau}$ converges smoothly to a \KE potential. In particular, 
$g^*_k \o_{\psi_{k\tau}}$ converges smoothly to a \KE metric.
\ethm

\begin{proof}
By Proposition \ref{weakProp} there exists $g_k \in G$ and a \KE potential 
$\psiKE \in \mathcal H_0$ such that $d_1(g_k.\psi'_k,\psiKE) \to 0$. We show below that in fact $g_k.\psi'_k \to_{C^\infty}\psiKE$. 

Focusing on the $\tau$-step Ricci iteration recursion, we can write:
\begin{flalign}
\big(g_{k+1}^{-1} \circ g_{k}\big)^* \Ric \ \o_{g_{k+1}.\psi'_{k+1}}& 
=
g_k^*\Ric \ \o_{\psi_{k+1}'}
=
g_k^* \Big(\frac{1}{\tau}\o_{\psi'_k} +\Big(1 - \frac{1}{\tau} \Big) \o_{\psi'_{k+1}}\Big) \nonumber\\
&= \frac{1}{\tau} \o_{g_k.\psi'_k} + \Big(1 - \frac{1}{\tau} \Big) \o_{g_k.\psi'_{k+1}}\nonumber\\
&= \frac{1}{\tau} \o_{g_k.\psi'_k} + \Big(1 - \frac{1}{\tau} \Big) \o_{(g^{-1}_{k+1} \circ g_k).g_{k+1}.\psi'_{k+1}}.\label{eq: pulled_it}
\end{flalign}
Set 
$$
\vp_k := g_k.\psi'_k \in \mathcal H_0
$$
 and 
$$
f_k:=g_{k}^{-1} \circ g_{k-1} \in G.
$$ 
With this notation, \eqref{eq: pulled_it} becomes:
\begin{flalign}
\Ric \ \o_{f_{k+1}. \vp_{k+1}}& = \frac{1}{\tau} \o_{\vp_k} + \Big(1 - \frac{1}{\tau} \Big) \o_{f_{k+1}.\vp_{k+1}}.\label{eq: pulled_it_new}
\end{flalign}

Without loss of generality we assume that $\o$ (the reference form) is 
\KEno. 
Using \eqref{eq: pulled_it_new} we can write:
$$
\i\ddbar\Big( \frac{1}{\tau}\vp_{k-1} + \Big(1-\frac{1}{\tau} \Big) f_k.\vp_{k} \Big) = \Ric \ \o_{f_{k}.\vp_{k}}
- \Ric \ \o = \i\ddbar \log \big ({\o^n}/{\o^n_{f_{k}.\vp_{k}}}\big).$$
This implies that 
$$\frac{1}{\tau}\vp_{k-1} + \Big(1-\frac{1}{\tau} \Big) f_k.\vp_{k} + \log (\o^n_{f_{k}.\vp_{k}}/ \o^n )=B_j \in \Bbb R.$$ 
Since $\log$ is a concave function, by Jensen's
inequality,
$$
\frac{1}{V}\int_M \log (\o^n_{f_{k}.\vp_{k}}/ \o^n ) \o^n \leq \log \frac{1}{V}\int_M  \o^n_{f_{k}.\vp_{k}} = 0.
$$
 By the triangle inequality, for $k$ sufficiently large, 
$$
d_1(0, \vp_{k-1})  \leq d_1(\psiKE,0) + 1.
$$ Using  \eqref{eq: d_1_double_est} we conclude that  $\int_M \vp_{k-1} \o^n \leq C$. These last two estimates combine to give 
$$B_j - \Big(1-\frac{1}{\tau} \Big) \frac{1}{V}\int_M f_k.\vp_{k}\o^n = \frac{1}{V}\int_M \vp_{k-1} \o^n + \frac{1}{V} \int_M \log (\o^n_{f_{k}.\vp_{k}}/ \o^n )\o^n \leq C.$$
Since $f_k.\vp_k \in \textup{PSH}(M,\o)$, it is well known that $\int_M f_k.\vp_{k}\o^n$ and $\sup_M f_k.\vp_{k}$ are comparable. As a result, 
$$B_j - \Big(1 - \frac{1}{\tau}\Big) \sup_M f_k.\vp_{k} \leq C,$$ 
hence we can write:
\begin{equation}\label{eq: main_MA_ineq}
\o_{f_k.\vp_k}^n = e^{B_j - (1-\frac{1}{\tau}) f_k.\vp_{k} -\frac{1}{\tau}\vp_{k-1} }\o^n \leq e^{C -\frac{1}{\tau}\vp_{k-1}}\o^n.
\end{equation}

Moreover, by Zeriahi's version of the Skoda integrability theorem \cite{Ze} (see \cite[Theorem 5.7]{DR} for a formulation that fits our context most), there exists $C>0$ such that, say,
$$
\intm e^{-\frac 3\tau\vp_{k-1}}\on\le C, \ k \in \Bbb N.
$$
Combining this estimate with \eqref{eq: main_MA_ineq}, we get that 
$$ 
||\o_{f_k.\vp_k}^n/\on||_{L^3(M,\on)}\le C.
$$
Now Ko\l odziej's estimate \cite{Bl,Ko} allows to conclude that 
the oscillation satisfies $\h{\rm osc}\, f_k.\vp_k\le C$ for some uniform $C$.
Note that for any $u\in\calH_0$,
it follows from \eqref{eq:AMest} that
$$
\inf u\le\frac{1}{V}\int u\o_u^n \le 0\le \frac{1}{V}\int u\o^n\le \sup u,
$$
so $u$ changes signs on $M$.
Thus, since $f_k.\vp_k\in\calH_0$, the oscillation bounds implies
a uniform bound
\beq
\lb{gkLinf2Eq}
||f_k.\vp_k||_{L^\infty(M)}\le C.
\eeq
Consequently, \eqref{eq: d_1_double_est}  yields
$$
d_1(0,f_k.\vp_k)=
d_1(f_k^{-1}.0,\vp_k)
\le C.
$$
Thus,
$$
d_1(f_k^{-1}.0,0)\le 
d_1(f_k^{-1}.0,\vp_k)+d_1(\vp_k,0)\le C'.
$$
By the arguments in the proof of \cite[Proposition 6.8]{DR}
(see also \cite[Lemma 2.7]{BDL} and \cite[Claim 7.11]{DR}), 
$\{f_k^{-1}\}_k$ is contained in a bounded set of $G$.
In particular, all derivatives up to order $m$, say, of $f_k^{-1}$
are bounded by some $C_m$ independently of $k$. So, 
to finish the proof, 
it suffices to estimate derivatives
of  
$$
h_k:=f_k.\vp_k
$$
(since that will imply the same estimates on
$f_k^{-1}. f_k. \vp_k=\vp_k$).

Note that $|\Delta_\o h_k|<C$ by the Chern--Lu argument of \cite[pp. 1539--1540]{R08}
since by \eqref{eq: pulled_it_new} we have
$$
\Ric \ \o_{h_{k+1}} = \Ric \ \o_{f_{k+1}.\vp_{k+1}} \geq \Big(1 - \frac{1}{\tau} \Big) \o_{f_{k+1}.\vp_{k+1}}= \Big(1 - \frac{1}{\tau} \Big)\o_{h_{k+1}}.
$$
 (cf. \cite[Corollary 7.8 (i)]{R14} with $C_1=0$ and $C_2=\big(\frac{1}{\tau} -1 \big)$).
The $C^{2,\alpha}$ and higher order estimates then follow the same way as in \cite{R08} (or by applying \cite[Theorem 5.1]{bl2} directly to \eqref{eq: main_MA_ineq}).

As we already have that $d_1(\varphi_k,\psiKE) \to 0$, an application of \eqref{eq: d_1_double_est} and the Arzel\`a-Ascoli compactness theorem finishes the argument.
\end{proof}

We note that in our arguments above the estimates depend on a positive lower bound to $\tau>0$. If this could be avoided, then one could hope that these estimates also hold in a scaled limit, as the iteration 
should converge to the K\"ahler--Ricci flow.


\paragraph{Acknowledgments.} Research supported by BSF grant 2012236, NSF grants DMS-1515703, DMS-1610202, and a Sloan Research Fellowship.

\vspace{0.1in}
\noindent {\sc University of Maryland}\\
{\tt tdarvas@umd.edu, yanir@umd.edu}
\end{document}